%% file: meas.tex
\def\**{$*\!*\!*$}
\input fibhead

\font\tenmsy=msym10

\textfont8=\tenmsy
\mathchardef\ssm="7872

\centerline{\** May 26, 1991, revised June 10, 1991 \**}\medskip

\centerline{\bf On the Lebesgue measure of the Julia set of
   a quadratic polynomial.}\medskip

\centerline{Misha Lyubich}\bigskip

\centerline{\bf \S1. Statement of the result.}\smallskip

The goal of this note is to prove the following

{\QP{\bf Theorem.} Let $p_a : z\mapsto z^2+a$ be a
quadratic polynomial which has no irrational indifferent periodic points,
and is not infinitely renormalizable.
Then the Lebesgue measure of the Julia set $J(p_a)$ is equal to zero.\medskip}

It was proven by McMullen that the Julia set of a cubic polynomial
$f_{a,b} : z\mapsto z^3-3a^2z+b$ has zero measure provided it is a Cantor
set (see [BH]). Our idea is to extract from the Julia set $J(p_a)$ an
essential part on which a ``renormalization" of $p_a$ is a
 ``polynomial-like map" $g$ (in a generalized sense). 
The construction involves the
Yoccoz partitions of the Julia set. Then using 
Branner-Hubbard-McMullen's method
and a notion of the modulus of a multiply connected domain, one can 
show that
the Julia set $J(g)$ has zero measure. By [L], $J(p_a)$ has zero measure as
well.

As part of the proof we discuss a property of the critical point
to be {\it persistently recurrent}, and relate our results to
corresponding ones for real one dimensional maps. 
In particular, we will show that
in the persistently recurrent case the restriction $p_a|\omega(0)$ is 
topologically minimal and has zero topological entropy. 

Let us mention that the Douady-Hubbard-Yoccoz rigidity theorem follows from 
the above result:

{\QP{\bf Corollary.} If $p_a$ is not infinitely renormalizable 
and has no attracting periodic points then it is  $J$-unstable.\medskip}
 
Indeed, by the Theorem, $p_a$ has no measurable invariant line fields on $J(p_a)$
and hence has no deformations concentrated on the Julia set [MSS]. It has no
deformations  on the Fatou set as well since there are no attracting
periodic points. So, $p_a$ is rigid. \smallskip

{\bf Acknowledgement.}
I would like to thank M. Shishikura for useful talks, after which we
both came  to the above theorem (but with different proofs).  I am 
grateful to J.~Milnor for helpful comments and discussions of the result.
I also take this opportunity to thank A. Douady and J. H. Hubbard
 for introducing me to the tableau technique.
\bigskip

\centerline{\bf \S2 Polynomial-like maps.}\smallskip
Let us introduce a notion of a {\it (generalized) polynomial-like map}
(compare [DH1]). Let $V$ and $V_i$, i=1,...,d, be open topological disks
with piecewise smooth boundaries such that the $V_i$ are
pairwise disjoint and cl$V_i\subset V$. 
(We don't require cl$V_i$ to be pairwise disjoint.)
A branched covering
$$g: \bigcup V_i\rightarrow V$$
will be called a (generalized) polynomial-like map. 

As usual, one can define the filled Julia set of $g$ as
$$K(g)=\{z: f^nz\in\cup V_i, n=0,1,...\},$$
and the Julia set $J(g)$ as its boundary.

In what follows we will assume that $g$ has a unique critical
point $c$ which is non-degenerate , and moreover that  $c\in K(g)$.

Let us consider  the nested sequence of inverse images $V^n=f^{-n}V$.
Clearly, cl$V^{n+1}\subset V^n$. Components $V^n_k$ are called $\it pieces$
of level $n$. Denote by $V^n(x)$ a piece of level $n$ containing the point
$x$. The pieces $V^n(c)$ will be called {\it critical}.

Following Branner and Hubbard [BH], for any $x\in K(g)$ one can organize 
the pieces in  a {\it tableau
} $T(x)=\{V^i(g^jx)\}_{i,j=0}^\infty$,
 mark the positions of the critical 
pieces and state three combinatorial rules for the resulting {\it marked grids}.

Now let us consider a topological disk $D$ containing
 a compact subset $K\subset D$. Let us
assign to the domain $A=D\ssm K$ the {\it modulus} $\mu(A)$ in the 
following way (see [A]). If $K$ has zero capacity then $\mu(A)=\infty$.
Otherwise
$$\mu(A)={1\over  I(u)}\eqno (1)$$
where $u$ is the harmonic function in $A$ which tends to 0 at regular points
 of $K$ and tends to 1 at regular points of $\partial D$ , and
$$I(u)=\int\int_A |{\rm grad}\; u|^2 dz d{\overline z}$$
is its Dirichlet integral .  Clearly, if we have  a $d$-sheeted branched covering
$h: (D,K)\rightarrow (D',K')$ then $\mu(A')=d\mu(A)$.

For a tableau $T(x)=V^n_k$ denote by $\mu^n_k$ the matrix of moduli
of the domains $A^n_k\equiv V^n_k\backslash V^{n+1}$. 

{\QP {\bf Lemma 1.}\item{ (i).} If a piece $V^n_k$ is not critical, $n>0,$ then
$\mu^{n-1}_{k+1}=\mu^n_k$. Otherwise $\mu^{n-1}_{k+1}=2\mu^n_k$.
\item{(ii).} If the critical tableau is aperiodic then
for any $x\in K(f)$ the nest of pieces $V^n(x)$ 
is of divergent type: 
$$\sum\mu(A^n(x))=\infty.$$\medskip}

{\bf Proof.} The first point is immediate from the properties of the moduli.
The second point is the formal consequence of the first one and the 
combinatorics of marked grids [BH].\QED

{\QP{\bf Corollary 1.} Assume that the critical tableau is
aperiodic.
Then $K(g)$ is a Cantor set.\medskip}

{\bf Proof.} Let us consider the annulus $A(x)=V\ssm\cap V^n(x) $
containing the disjoint union of domains $A^n(x).$ Since the embeddings
$A^n(x)\subset A(x)$ are homotopically non-trivial, the Gr\"{o}tzsch inequality
yields
$$\mu (A(x))\geq \sum\mu(A^n(x))=\infty.$$ 
So, $\cap V^n(x)$ is a point.\QED

Now let us state an analytical lemma which generalizes McMullen's one
onto multiply connected domains (see [BH], \S 5.4).
Denote by $\lambda$ the Lebesgue measure on the plane.\smallskip

{\QP{\bf Lemma 2.} Let $D$ be a topological disk and $K\subset D$ be its
compact subset consisting of finitely many components, $A=D\ssm K.$ Then
    $${\lambda(D)\over \lambda(K)}\geq 1+4\pi\mu(A).$$\medskip}
{\bf Proof.} The modulus $\mu(A)$ can also be defined as the reciprocal to the
extremal length $\sigma(\Gamma)$ of the family $\Gamma$ of (non-connected)
curves separating $K$ from $\partial D$. For such a curve $\gamma$ there is a
decomposition of $K$ into the union of disjoint pieces $K_i$ 
 surrounded by pieces $\gamma_i$ of $\gamma.$ Then
$$\sigma(\Gamma)\geq\inf_{\gamma\in\Gamma}{|\gamma|^2\over\lambda(A) }\geq
\inf_{\gamma\in\Gamma}{\sum |\gamma_i|^2\over \lambda(A)}\geq
4\pi{\sum\lambda(K_i)\over \lambda(A)}=4\pi{\lambda(K)\over\lambda(A)}$$
where the last inequality follows from the isoperimetric one. Now the
required estimate follows.\QED

{\QP{\bf Corollary 2.} If the critical tableau is aperiodic then the Lebesgue
measure of $K(g)$ is equal to zero.\medskip}

{\bf Proof.} One can repeat the McMullen's argument word by word. Just for
fun we will slightly modify it. 

Let us organize the set of pieces $V^n_k$ in a tree joining $V^n_k$ with
$V^{n+1}_i$ in the case when $V^{n+1}_i\subset V^{n}_k$. Let us assign to each
edge $[U,W]$ of the tree a number
$$\nu[U,W]\equiv\nu(U)=\min(\mu(U),1/2),$$
and to each branch $\gamma$ a number $\nu(\gamma)$ which is the sum of
 $\nu[U,W]$
over all edges of $\gamma$. Denote by $\Gamma_n$ the family of all branches
of length $n,\; n\leq\infty$
(saying ``branch" we mean a path in the tree beginning at the root vertex $V$).
By Lemma 1, $$\nu(\gamma)=\infty\eqno (2)$$ for any $\gamma\in \Gamma_\infty.$
 Let us show that 
$$M_n\equiv\min_{\gamma\in\Gamma_n}\nu(\gamma)\to\infty.\eqno (3)$$

Indeed, given a $C$, consider a subtree of vertices  $W$ such that 
 $\nu[V,W]\leq C$ where $[V,W]$ is the branch ending at $W$.
 Now use so called
K\"{o}nig lemma: if a tree with finitely many branches at any vertex has
arbitrary long branches then it has an infinite branch. Along this branch the
divergent (2) condition fails.

 By Lemma 2, for any vertex $U$ of level $n$
$${\lambda(V^{n+1}\cap U)\over\lambda(U)}
\leq\exp(-b\nu(U))$$
with an appropriate constant $b$.
Now one can easily derive from here  ( induction in $n$ ) that 
$$\lambda(V^n)\leq\exp(-b M_n)\lambda(V),$$
and  by (3) this goes down to 0. \QED
\bigskip

\centerline{\S3 \bf Persistent recurrence and renormalization. }
\medskip

Set $c=0$, the critical point of a quadratic polynomial $f=p_a:z\to z^2+a$.
If $f$ has an attractive periodic point
 or a rational
indifferent periodic point then $\lambda(J(f))=0$ ([DH2], [L]). 
So, we will  assume till the end of the paper
that $f$ has no such points (and is not infinitely renormalizable). 
Then we are exactly in the situation
studied recently by Yoccoz. We will use the Yoccoz construction without
detailed explanation (see [H]).

For (open) pieces $V^n_k$ of the Yoccoz partitions
we will use the same notations as for the
Branner-Hubbard ones. By $\partial V^n$ denote the union of $\partial V^n_k$.
The intersection $\partial V^n\cap J(f)$ consists of finitely
many points, preimages of a fixed point $\alpha$. Set 

$$V^n(z)={\rm int}\bigcup{\rm cl}V^n_k$$ 
over the pieces $V^n_k$ of level $n$ whose closure contain $z$ (if $z$ is not
a preimage of $\alpha$, $V^n(z)$ is just the piece of level $n$ containing $z$).
 The  $V^n(z)$ is a neighborhood of $z$. It follows from  the Yoccoz Theorem
that
$$ {\rm diam}\;V^n(z)\to 0\qquad {\rm as}\; n\to\infty \eqno (4)$$
uniformly in $z$.

 Given a domain $U$  and an orbit $\{z_k\}_{k=0}^n$ such that
$z_n\in U$, one can {\it pull } $U$ {\it back} along this orbit, 
that is to consider
the string of domains $U_k,\; k=0,1,...,n,$ such that $U_n=U$, and 
$U_{k}$ is the component of $f^{-1} U_{k+1}$ containing $z_k$. 
In particular, if $U=V^l(z_n)$ 
then  $U_k=V^{l+n-k}(z_k)$ . 

The {\it order} ord${\bf U}$ of the pull-back ${\bf U}$
is the number of domains
$U_k$ containing the critical point $c$. 
The pull-back of zero order (that is,  none of $U_k$  covers the 
critical point) will be called {\it univalent} .
 
{\QP{\bf Lemma 3.} Let $W$ and $U$ be two domains intersecting the Julia set
$J(f)$, and $f^n(W)\subset U$ for some $n\in {\bf N}$. Then
$${\rm diam}\;U<\delta \Rightarrow {\rm diam}\;W<\epsilon_n(\delta)<
\epsilon(\delta)$$
with $\epsilon (\delta)\to 0$ as $\delta\to 0$,
and $\epsilon_n(\delta)\to 0$ as $n\to\infty$.\medskip}

{\bf Proof.} Fix $\epsilon>0$. According to (4), there is an $\ell\in{\bf N}$
such that
$${\rm diam}\;V^\ell(z)<\epsilon\eqno (5)$$
for any $z\in J(f)$. The sets $V^\ell(z)$ form an open covering of the 
Julia set. Let $\delta $ be its {\it Lebesgue number}. This means that any
set $U$ of diameter less than $\delta$ is contained in some domain 
$V^\ell(\zeta),\qquad \zeta\in J(f).$

Let us apply this to a given domain $U$, and find $z\in W\cap J(f)$ such that
$f^nz=\zeta$. Pulling $V^\ell$ back along the orbit of $z$ we come to a 
piece $V^{\ell+n}(z)$ containing $W$. 
Now  (5) yields that ${\rm diam}\;W<\epsilon$,
and we are done.~~\QED

Denote by $B(z,r)$ the Eucledian disk of radius $r$ centered at $z$.\smallskip

{\QP{\bf Lemma 4.} If $\lambda(J(p_c))>0 $ 
then for almost all $z\in J(p_c)$ we have 
$$\omega(z)=\omega(c)\ni c.$$
So, $\lambda(J(p_c))=0$ provided $c$ is non-recurrent.\medskip}

{\bf Proof.} The inclusion $\omega(z)\subset\omega(c)$ for almost all $z$ 
follows from the Koebe Distortion Theorem (see [L]). We are going to
 show that $\omega(z)\supset c$ for almost all $z$ which certainly
 implies the required statement. Let
$${\rm dist}({\rm orb(z)},c)\geq\gamma. \eqno (6)$$
Find a $\delta>0$ such that $\epsilon (\delta)<\gamma$ (see Lemma 3).
Consider the disk $B(z_n,\delta)$ around $z_n\equiv f^nz$, and let 
${\bf U}=\{U_k\}_{k=0}^n$ be its pull-back
along the orbit of $z$. By Lemma 3, diam$U_k<\gamma,\; k=0,...,n$,
hence the pull-back ${\bf U}$ is univalent.
 By the Koebe Theorem, $z$ is not a density
point of $J(p_c)$, and the statement follows. \QED

{\bf Remark.} Consider a set $J_\gamma$ of all $z$ satisfying (6).
Then the restriction $f|J_\gamma$ is {\it expanding}: there exist $C>0$
and $q>1$ such that
$$(f^n)'(z)\geq cq^n,\qquad z\in J_\gamma,\; n=0,1,...$$
Indeed, by Lemma 3, diam$U_0<\epsilon_n(\delta)\to 0$ as $n\to\infty$ 
(where $U_0$ is the set from  the above proof). Hence, there is an $n$ 
such that $(f^n)'(z)\geq q>1$ for all $z\in J_\gamma$,
and the statement follows. The analogous
fact is well-known in one-dimensional dynamics [G]. 

Let us consider now the Yoccoz $\tau$-function. For an $n\in {\bf N}$ it
assigns the biggest $m\in [0,n-1]$ for which the piece $f^{n-m} V^n(c)$
of level $m$ contains
the critical point $c$ (if there is no such an $m$, set $\tau(n)=-1$.)
The critical point is called {\it persistently recurrent} if 
$\tau(n)\to\infty$ as $n\to\infty.$ The following lemma was probably
 known to several people.
 
{\QP{\bf Lemma 5.} If the critical point is not persistently recurrent then
$\lambda(J(p_c))=0.$ \medskip} 

{\bf Proof.} Since $c$ is not persistently recurrent, there is an $N$ and
arbitrary large $l$ such that $f^l$ is a double covering of $Q\equiv V^{N+l}(c)$
over $V^N(c).$ Let
$$z\in J(f)\ssm \bigcup_{n=0}^\infty f^{-n}\alpha,$$
and consider the first moment $n$ for which $f^nz\in Q.$ 
Let  ${\bf Q}=Q_0,...,Q_n$ be the pull-back of $Q$ along the orbit of $z$.
Since the level of $Q_i$ is bigger than $N+l$ for $i<n$, these $Q_i$ don't
cover $c$ which means that the pull-back is univalent. 

The boundary $\partial V^N$ consists of $N$-preimages of invariant rays through
$\alpha$ and $N$-preimages of an equipotential level.
It follows that there exists a $\delta>0$ such that for any $u\in V^N(c)$ 
there is a $k$ and a neighborhood $U\subset V^N(c)$ around $u$ such that
$f^k$ univalently maps $U$ onto a disk $B(f^ku,\delta) $of radius $\delta$ 
centered at $f^ku$.

 Set $u=f^{n+l}z$ and find corresponding $U$ and $k$. Now consider two cases

{\item (i)} $f^l c$ does not belong to $U$ or $|f^{l+k}c-f^ku|>\delta/100$. 
Then the pull-back $T_i$ of the
disk $B(f^ku,\delta/200)$ along the orbit $\{z_m\}_{m=0}^{n+l+k}$ is univalent.
By the Koebe theorem,  the density of the Julia set in $T_0$ is bounded away
from 1. Since $l$ is arbitrary large, $z$ is not a density point of $J(p_c)$. 

{\item (ii)} $f^l c\in U $ and $|f^{l+k}c-f^ku|<\delta/100$. 
Then we can find a disk $B_1$  centered
at the critical value and such that its $f^{l+k-1}$-image 
lies in between $B(f^k(u),\delta/10)$ and $B(f^k(u),\delta/2)$.
By the Koebe Theorem, the density of $J(p_c)$ at $B_1$ is bounded away from 1.
Now let us consider the preimage $B_0$ of $B_1$. It is a disk centered at
$c$ and contained in $V^{N+l}(c)$. The squaring map in the disk cannot expand 
density of thin sets too much. It follows that the density of the Julia set
in $B_0$ is bounded away from 1 as well. Now pulling a little bit smaller
disk along the orbit $\{z_i\}_{i=0}^n$, 
we conclude again that $z$ is not a density
point of the Julia set.  \QED 

Now let us concentrate on the main case of a persistently recurrent critical
point. 

{\QP{\bf Lemma 6.} The following properties are equivalent:
{\item (i).} The critical point $c$ is persistently recurrent.
{\item (ii)}.{\it Absence of long univalent pull-backs.}
 Given an $\epsilon>0$, there is an $N\in {\bf N}$ with the
following property. Let $\overline{z}=\{z,z_{-1},...,z_{-n}\}$ be any backward
orbit in $\omega(c)$ . If the pull-back of the disk $B(z,\epsilon)$ along
$\overline{z}$ is univalent then $n\leq N$.\medskip}

{\bf Remark.} Property (ii) equivalent to persistent recurrence appeared in
one dimensional setting in [BL]. Note that the formulation does not involve
any particular partitions of the Julia set.\medskip 

{\bf Proof.} (i)$\Rightarrow$ (ii). Find an $\ell$ so that (5) holds.
Take a backward orbit $\overline{z}  $ and assume that the pull-back
of $B(z,\epsilon)$ along it is univalent. By (5), 
$B(z,\epsilon)\supset V^\ell(z)$, hence the pull-back of $V^\ell(z)$ along
$\overline{z}$ is also univalent. 

 Since $z_{-n}\in \omega(c)$, there is
an $s\in {\bf N}$ for which 
$$c_s\in V^{\ell+n}(z_{n})\eqno (7).$$ 
So, we can pull the piece $V^{\ell+n}(z_{-n})$ 
back along the orbit of $c$ till the first moment $t$
when it covers the critical
point (if $s$ is the first moment for which (7) holds then $t=s$).
We obtain a critical piece $V^{\ell+n+t}(c)$ such that
$$f^{n+t}:V^{\ell+n+t}(c)\rightarrow V^\ell$$
 is a double covering. By persistent recurrence, we get a uniform bound on $n$.
\smallskip
(ii)$\Rightarrow$ (i). If $f$ is not persistently recurrent, we can find
a critical piece $V^\ell(c)$ allowing arbitrarily long univalent pull-backs
along $\omega(c)$. Since $V^\ell(c)\supset B(c,\epsilon)$ for some $\epsilon>0$,
we arrive at a contradiction.\QED
\medskip 
The following Corollary seems to be interesting by itself. It will not be used
for the proof of the Theorem.\smallskip

{\bf Corollary}(cf [BL],\S 11). If $f$ is persistently recurrent then
{\item (i).} The restriction $f|\omega(c)$ is topologically minimal
 (that is, all orbits
are dense in $\omega(c)$)
{\item (ii).} $f|\omega(c)$  has zero topological entropy.
\medskip
{\bf Proof.} (i).Let $z\in\omega(c)$. We should prove that $c\in\omega(z)$.
Otherwise dist(orb$(z),c)\geq\gamma>0$. By Lemma 3, for sufficiently small
$\delta$, the pull-back of the disk $B(z_n,\delta)$ along the orbit of $z$
is unimodal, contradicting lemma 6(ii).\smallskip
(ii) If there is a measure $\mu$ of positive entropy supported on $\omega(c)$,
then the Pesin theory of unstable manifolds yields the existence of a disk
$B(z,\delta),\;z\in J(f)$, having an infinite pull-back along $\omega(c)$
contradicting Lemma 6(ii) again.\QED

We are prepared to prove the main lemma.

{\QP{\bf Lemma 7.} Let $c$ be persistently recurrent.
Then there is a polynomial-like map
$$g:\bigcup V_i\rightarrow V$$
such that 
{\item (i)} $g|V_i=p_c^{n_i}$ for some $n_i$.
{\item (ii)} $c$ is the unique critical point of $g$.
{\item (iii)} $c\in K(g).$ \medskip}

{\bf Proof.}  Let us take a non-degenerate annulus
$$V^{n-1}(c)\ssm V^n(c)$$
around the critical point (see [H]), so that cl$V^n(c)\subset V^{n-1}(c)$.
Set $V=V^n(c)$ Then
$$\partial(f^jV)\cap {\rm cl}V=\emptyset,\qquad j=1,2,...\eqno (8).$$
Consider all returns $c_{m(i)}$ of orb$c$ into $V$, and let $l(i)=m(i+1)-m(i).$
Consider all pull-backs $\{f^k V_i\}_{k=0}^{l(i)}$ of $V$ along the strings 
$\{c_k\}_{k=m(i)}^{m(i+1)}.$ For any $i$ all intermediate pieces 
$f^kV_i, \; 0<k<l(i) ,$ lie outside $V$, so the maps 
$f^{l(i)-1}: fV_i\rightarrow V$ are univalent.  Hence the maps
$f^{l(i)}: V_i\rightarrow V$ are either univalent (if the piece $V_i$
 is not critical) or double coverings.

By Lemma 6(ii), $l_i$ are uniformly bounded: $l_i\leq L.$
So, actually we have only finitely many different sets $V_i.$
Let us show that these sets are disjoint. Indeed, otherwise $V_i\supset V_j$
for two different sets $V_i$ and $V_j$. Let us push $V_i$ forward till the
level $n$:  $f^k V_i=V. $ Then $f^k V_j\subset V$ despite the fact that all
sets $f^k V_j$ of level greater than $n$ lie outside of $V$.
(Actually, (8) yields more: cl$V_i\cap$cl$V_j=\emptyset$ for $i\not=j$.)

Now let us show that cl$V_i\subset V$. Indeed, otherwise
$\partial V_i\cap \partial V\not=\emptyset$. Let $\ell>n$ be the level of
$V_i$, $j=\ell-n$. Then 
$$\partial V\cap\partial(f^jV)\supset
f^j(\partial V_i\cap\partial V)\not=\emptyset$$
contradicting (8).

We have shown that $g|\bigcup V_i$ is a generalized polynomial-like map.
Since $g$ is non-univalent only on the critical piece, $c$ is the only
critical point of $g$. Since \break $gc_{m(i)}=c_{m(i+1)}\in V_{i+1}$,
$c\in K(g)$. \QED

{\bf Proof of the Theorem.} If $f$ is finitely renormalizable, let us 
renormalize it to be non-renormalizable. It follows from Lemma 4 that the
measure of the Julia set of the original polynomial and its renormalization
is simultaneously positive or zero. So, we will assume in what follows that
$f$ is non-renormalizable (and $c$ is persistently recurrent).

Consider the polynomial-like map $g$ constructed in the previous lemma. 
Since $p_c$ is non-renormalizable, $g$ has an aperiodic tableau.
By Corollaries  1 and 2, $K(g)$ is a Cantor set of zero measure.
Let $z$ be a typical point of $J(f)$, so that the conclusion of Lemma 4 holds.
Consider the  moments $n(i)$ when the orb($z$) returns  to $V$.
Since $\omega(c)\cap V\subset \cup V_i$,  eventually $z_{n(i)}\in \cup V_i$.
Hence $f^n z\in K(g)$ for some $n$. Now the Theorem follows.\QED
\bigskip

\centerline{\bf References.}\medskip

{\item [A]} L. V. Ahlfors. Conformal invariants. McGraw Book Company, 1973.
\smallskip
{\item [BH]} B. Branner and J. H. Hubbard. The iteration of cubic polynomials.
  Part II. Matematik Institut, Denmark, 1989-12.\smallskip
{\item[BL]} A. Blokh and M. Lyubich. Measurable dynamics of S-unimodal
maps of the interval. Preprint, SUNY, Stony Brook, 1990/2.
{\item [DH1]}  A. Douady and J. H. Hubbard. On the dynamics of polynomial-like
   mappings. Ann. Sci. l'\'{E}cole Norm. Sup., 18 (1985), 287-343.\smallskip
{\item [DH2]}  A. Douady and J. H. Hubbard. 
\'{E}tude dynamique des polyn\^{o}mes complexes. Part I.
  Preprint Orsay,  84-02. \smallskip
{\item [G]} J. Guckenheimer. Sensitive dependence to initial conditions
for one-dimensional maps. Comm. Math. Phys., 70 (1979), 133-160. 
{\item [H]} J. H.Hubbard, according J.-C. Yoccoz. Puzzles and quadratic tableaux.
  Manu\-script.\smallskip
{\item [L]} M. Lyubich. On typical behavior of the trajectories of a rational
   mapping of the sphere. Soviet Math. Dokl., 27 (1983), n 1, 22-25.\smallskip
{\item [MSS]} R. Ma\~{n}\'{e}, P. Sad and D. Sullivan. On the dynamics of
  rational maps. Ann. Sci. l'\'{E}cole Norm. Sup., 16 (1983), 193-217.

\end

%% file: fibhead.tex

 at 12truept
\font\tenmsy=msym10
\textfont8=\tenmsy
\mathchardef\ssm="7872

\input psfig
\mathsurround = 2pt
\abovedisplayskip=6pt
\belowdisplayskip=6pt

\def \QP{\narrower\medskip\noindent}
\def \QED {\rlap{$\sqcup$}$\sqcap$\smallskip}

\def\ref{\hangindent=1pc \hangafter=1 \noindent}